\documentclass[12pt,reqno]{amsart}

\usepackage{amssymb}
\usepackage{bm}
\usepackage{fixmath}
\usepackage{mathrsfs}
\usepackage{color}

\newtheorem{theorem}{Theorem}[section]

\newtheorem{lemma}[theorem]{Lemma}
\newtheorem{claim}[theorem]{Claim}

\numberwithin{equation}{section}

\begin{document}

\title[]{On algebraic differential equations concerning the Riemann-zeta function and the Euler-gamma function}

\author[]{Qiongyan Wang$^1$\quad Manli Liu$^2$*\quad Nan Li$^3$}
\address{$^1$ School of Mathematical sciences, Peking University, Beijing, 100871, P.R. China
\vskip 2pt $^2$ School of Mathematics, Shandong University, Jinan, Shandong, 250100, P.R. China
\vskip 2pt $^3$ School of Mathematics, Qilu Normal University, Jinan, Shandong, 250013, P. R. China
\vskip 2pt \hspace{1.5mm} Email:{\sf qiongyanwang@math.pku.edu.cn, lml6641@163.com, nanli32787310@163.com}}

\thanks{{\sf Manli Liu* is the corresponding author.}}
\thanks{The research was partially supported by NSFC of China(No. 11801215), NSFC of Shandong (No. ZR2018MA014).}
\begin{abstract}
In this paper, we prove that $\zeta$ is not a solution of any non-trivial algebraic differential equation whose coefficients are polynomials in $\Gamma, \Gamma^{(n)}, \Gamma^{(l)}$ over the ring of polynomials in $\mathbb{C}$, $l>n\geq 1$ are positive integers. We extended the result that $\zeta$ does not satisfy any non-trivial algebraic differential equation whose coefficients are polynomials in $\Gamma, \Gamma', \Gamma''$ over the field of complex numbers, which is proved by Li and Ye\cite{8}.
\end{abstract}

\maketitle

\section{Introduction and main results}
It is known \cite{6} that the Euler gamma-function
$$\Gamma(z)=\int_{0}^{+\infty}t^{z-1}e^{-t}dt, ~~~ \Re z> 0$$
(analytically continued as a meromorphic function in the complex plane $\mathbb{C}$) does not satisfy any non-trivial algebraic differential equation whose coefficients are polynomials in $\mathbb{C}$. That is, if $P(v_{0}, v_{n}, ..., v_{l})$ is a polynomial of $n+1$ variables with polynomial coefficients in $z\in\mathbb{C}$ such that
$P(\Gamma, \Gamma', ..., \Gamma^{(n)})(z)\equiv 0$ for $z\in\mathbb{C}$, then $P\equiv 0$. Hilbert\cite{5}, in Problem 18 of his famous list of 23 problems, stated that Riemann zeta function
$$\zeta(z)=\sum_{n=1}^{+\infty}\frac{1}{n^{z}}$$
does not satisfy any non-trivial algebraic differential equation whose coefficients are polynomials in $\mathbb{C}$, the problem was solved in great generality by Ostrowski\cite{11}. Since then, these results have been extended in many different directions(see e.g. \cite{2}, \cite{3}, \cite{7}, \cite{12}), to list a few. It is well known that $\zeta$ is associated with $\Gamma$ by the Riemann functional equation
\begin{equation}\label{eq0.1}
\begin{aligned}
\zeta(1-z)=2^{1-z}\pi^{-z}\cos \frac{\pi z}{2}\Gamma(z)\zeta(z).
\end{aligned}
\end{equation}

By virtue of (\ref{eq0.1}), one knows from Bank-Kaufman \cite{1} and Liao-Yang \cite{9} that neither $\zeta$ nor $\Gamma$ satisfy any non-trivial algebraic differential equation whose coefficients are meromorphic functions, growing strictly slower than $e^{z}$ as $|z|\rightarrow +\infty$ or having period $1$, detailed description of this matter can be found in Li and Ye \cite{8}. The result \cite{1} can not applied to algebraic differential equations involving both $\zeta$ and $\Gamma$ simultaneously, since $\zeta$ and $\Gamma$ grow strictly faster than $e^{z}$ as $|z|\rightarrow +\infty$; see Ye \cite{16}.

Recently, Markus \cite{10} proved that $\zeta(\sin(2\pi z))$ cannot satisfy any non-trivial algebraic differential equations whose coefficients are polynomials in $\Gamma$ and its derivatives, and he conjectured that $\zeta$ itself cannot satisfy any non-trivial algebraic differential equations whose coefficients are polynomials in $\Gamma$ and its derivatives. Thus, we are interested in knowing whether there is a non-trivial polynomial $P(u_{0}, u_{1}, ..., u_{m}; v_{0}, v_{1}, ..., v_{n})$
such that
$$P(\zeta, \zeta', ..., \zeta^{(m)}; \Gamma, \Gamma', ..., \Gamma^{(n)})(z)\equiv 0, ~~~ z\in\mathbb{C}.$$
In other words, one wants to know whether $\zeta$ satisfies any non-trivial algebraic differential equation whose coefficients are differential polynomials of $\Gamma$.

In this paper, we prove the main result as follows.
\begin{theorem}\label{M-thm}
Let $m, n, l$ be positive integers, $n<l$. $P(u_{0}, u_{1}, ..., u_{m}; v_{0}, v_{n}, v_{l})$ is a polynomial of $m+4$ variables with polynomial coefficients in $z\in\mathbb{C}$ such that
\begin{equation}\label{eq1.1}
\begin{aligned}
P(\zeta, \zeta', ..., \zeta^{(m)}; \Gamma, \Gamma^{(n)}, \Gamma^{(l)})(z)\equiv 0
\end{aligned}
\end{equation}
for $z\in\mathbb{C}$. Then the polynomial $P$ must be identically equal to zero.
\end{theorem}
Theorem~\ref{M-thm} shows that $\zeta$ does not satisfy any non-trivial algebraic differential equation whose coefficients could be polynomials in $\Gamma, \Gamma^{(n)}, \Gamma^{(l)}$. Our theorem extended the result that $\zeta$ does not satisfy any non-trivial algebraic differential equation whose coefficients could be polynomials in $\Gamma, \Gamma', \Gamma''$, which is proved by Li and Ye \cite{8}. As of \cite{8B},using its terminology, only $\zeta$ and the family of distinguished polynomials of $\Gamma, \Gamma', ..., \Gamma^{(n)}$ were considered.

\section{Preliminary}
To prove our theorem, we need the following celebrated theorem from Voronin \cite{14}.
\begin{lemma}\label{lem2.1}
	Fix $x\in\left(\frac{1}{2}, 1\right)$ for $z=x+i y\in\mathbb{C}$. Define
$$\gamma(y):=(\zeta(x+i y), \zeta'(x+i y), ..., \zeta^{(m)}(x+i y))$$
to be a curve in $y$. Then, $\gamma(\mathbb{R})$ is everywhere dense in $\mathbb{C}^{m+1}$.
\end{lemma}

\section{Proof of Theorem~\ref{M-thm}}
Let
$P(u_{0}, u_{1}, ..., u_{m}; v_{0}, v_{n}, v_{l})$
be a polynomial in its arguments with coefficients in the field $\mathbb{C}$ such that (\ref{eq1.1}) is satisfied. In view of the discussions in Section $2$ of \cite{10}, one can simply assume that the coefficients of $P$ are constants. Denote by
\begin{equation*}
\begin{aligned}
\Lambda:=\{\lambda:=(\lambda_{0}, \lambda_{n}, \lambda_{l}): \lambda_{0}, \lambda_{n}, \lambda_{l} ~\text{are non-negative integers}~\},
\end{aligned}
\end{equation*}
a triple-index set having a finite cardinality, and define
\begin{equation*}
\begin{aligned}
\Lambda_p:=\{\lambda\in\Lambda:|\lambda|=p ~\text{with}~ |\lambda|:=\lambda_{0}+\lambda_{n}+\lambda_{l}\},
\end{aligned}
\end{equation*}
\begin{equation*}
\begin{aligned}
\Lambda_{q}^\ast:=\{\lambda\in\Lambda:|\lambda|^\ast=q ~\text{with}~ |\lambda|^\ast:=n\lambda_{n}+l\lambda_{l}\}.
\end{aligned}
\end{equation*}
Then, there is a non-negative integer $L$ such that
$$P(u_{0}, u_{1}, ..., u_{m}; v_{0}, v_{n}, v_{l})=\sum_{p=0}^{L}\sum_{\lambda\in\Lambda_{p}}a_{\lambda}(u_{0}, u_{1}, ..., u_{m})v_{0}^{\lambda_{0}}v_{n}^{\lambda_{n}}v_{l}^{\lambda_{l}},$$
where $a_{\lambda}(u_{0}, u_{1}, ..., u_{m})$ is a polynomial of $m+1$ variables with coefficients in $\mathbb{C}$. Set, for each $p=0, 1, ..., L$,
$$P_{p}(u_{0}, u_{1}, ..., u_{m}; v_{0}, v_{n}, v_{l})=\sum_{\lambda\in\Lambda_{p}}a_{\lambda}(u_{0}, u_{1}, ..., u_{m})v_{0}^{\lambda_{0}}v_{n}^{\lambda_{n}}v_{l}^{\lambda_{l}}.$$
Arranging $P_{p}(u_{0}, u_{1}, ..., u_{m}; v_{0}, v_{n}, v_{l})$ in $v:=(v_{0}, v_{n}, v_{l})$ in the ascending order of $q=|\lambda|^{\ast}$, we get that there is a non-negative integer $M_{p}$ such that
\begin{equation}\label{eq1.2}
\begin{aligned}
P_{p}(u_{0}, u_{1}, ..., u_{m}; v_{0}, v_{n}, v_{l})=\sum_{q=0}^{M_{p}}\sum_{\lambda\in\Lambda_{p}\cap\Lambda_{q}^{\ast}}a_{\lambda}(u_{0}, u_{1}, ..., u_{m})v_{0}^{\lambda_{0}}v_{n}^{\lambda_{n}}v_{l}^{\lambda_{l}}.
\end{aligned}
\end{equation}
Consequently,
\begin{equation}\label{eq1.3}
\begin{aligned}
P(u_{0}, u_{1}, ..., u_{m}; v_{0}, v_{n}, v_{l})=\sum_{p=0}^{L}\sum_{q=0}^{M_{p}}\sum_{\lambda\in\Lambda_{p}\cap\Lambda_{q}^{\ast}}a_{\lambda}(u_{0}, u_{1}, ..., u_{m})v_{0}^{\lambda_{0}}v_{n}^{\lambda_{n}}v_{l}^{\lambda_{l}}.
\end{aligned}
\end{equation}
\begin{claim}\label{cla 3.1}
Assume (\ref{eq1.1}) holds, then, for each $0\leq p\leq L$ and $z\in \mathbb{C}$, one has
\begin{equation*}
P_{p}(\zeta, \zeta', ..., \zeta^{(m)}; \Gamma, \Gamma^{(n)}, \Gamma^{(l)})(z)\equiv 0.
\end{equation*}
\end{claim}
\emph{Proof.}
Suppose $p_0$ is the smallest index among $\{0, 1, ..., L\}$ such that
\begin{equation*}
P_{p_0}(\zeta, \zeta', ..., \zeta^{(m)}; \Gamma, \Gamma^{(n)}, \Gamma^{(l)})(z)\not\equiv 0.
\end{equation*}
Then, we have
\begin{equation}\label{eq1.4}
\begin{aligned}
P_{p_0}\left(\zeta, \zeta', ..., \zeta^{(m)}; 1, \frac{\Gamma^{(n)}}{\Gamma}, \frac{\Gamma^{(l)}}{\Gamma}\right)(z)=\frac{P_{p_0}(\zeta, \zeta', ..., \zeta^{(m)}; \Gamma, \Gamma^{(n)}, \Gamma^{(l)})(z)}{\Gamma^{p_0}(z)}\not\equiv 0.
\end{aligned}
\end{equation}
Define the digamma function $f:=\frac{\Gamma'}{\Gamma}$, and introduce inductively
\begin{equation*}
\begin{aligned}
\frac{\Gamma'}{\Gamma}=&f\\
                      =&f\big[1+\frac{f'}{f^2}(c_1+\varepsilon_1)\big] ~ \text{for}~ c_1=0 ~ \text{and} ~ \varepsilon_1=0,\\
\frac{\Gamma''}{\Gamma}=&\left(\frac{\Gamma'}{\Gamma}\right)'+\left(\frac{\Gamma'}{\Gamma}\right)^{2}=f'+f^{2}\\
                       =&f^{2}\big[1+\frac{f'}{f^2}(c_2+\varepsilon_2)\big] ~ \text{for}~ c_2=1 ~ \text{and} ~ \varepsilon_2=0,\\
\frac{\Gamma'''}{\Gamma}=&\left(\frac{\Gamma''}{\Gamma}\right)'+\frac{\Gamma''}{\Gamma}\cdot\frac{\Gamma'}{\Gamma}=f''+3ff'+f^{3}\\
                        =&f^{3}\big[1+\frac{f'}{f^2}(c_3+\varepsilon_3)\big] ~ \text{for}~ c_3=3 ~ \text{and} ~ \varepsilon_3=\frac{f''}{ff'},\\
\frac{\Gamma^{(4)}}{\Gamma}=&\left(\frac{\Gamma'''}{\Gamma}\right)'+\frac{\Gamma'''}{\Gamma}\cdot\frac{\Gamma'}{\Gamma}=f'''+4ff''+3(f')^{2}+6f^{2}f'+f^{4}\\
                           =&f^{4}\big[1+\frac{f'}{f^2}(c_4+\varepsilon_4)\big] ~ \text{for}~ c_4=6 ~ \text{and} ~ \varepsilon_4=\frac{f'''}{f^{2}f'}+4\frac{f''}{ff'}+3\frac{f'}{f^{2}},\\
\frac{\Gamma^{(5)}}{\Gamma}=&\left(\frac{\Gamma^{(4)}}{\Gamma}\right)'+\frac{\Gamma^{(4)}}{\Gamma}\cdot\frac{\Gamma'}{\Gamma}\\
                           =&f^{4}+5ff'''+10f'f''+10f^{2}f''+15f(f')^{2}+10f^{3}f'+f^{5}\\
                           =&f^{5}\big[1+\frac{f'}{f^2}(c_5+\varepsilon_5)\big] ~ \text{for}~ c_5=10 ~ \text{and} \\
                               &\varepsilon_5=\frac{f^{(4)}}{f^{3}f'}+5\frac{f'''}{f^{2}f'}+10\frac{f''}{ff'}+10\frac{f''}{f^{3}}+15\frac{f'}{f^{2}},\\
\vdots\\
\frac{\Gamma^{(n)}}{\Gamma}=&f^{n}\big[1+\frac{f'}{f^2}(c_n+\varepsilon_n)\big]~\text{upon assumption, and then we can deduce that}\\
\frac{\Gamma^{(n+1)}}{\Gamma}=&\left(\frac{\Gamma^{(n)}}{\Gamma}\right)'+\frac{\Gamma^{(n)}}{\Gamma}\cdot\frac{\Gamma'}{\Gamma}
                             =f^{n+1}+f^{n-1}f'(c_n+n+\varepsilon_n)\\
                              &+f^{n-2}f'\varepsilon'_n +[(n-2)f^{n-3}(f')^{2}+f^{n-2}f''](c_n+\varepsilon_n)\\
                             =&f^{n+1}\big[1+\frac{f'}{f^2}(c_{n+1}+\varepsilon_{n+1})\big] ~ \text{for}~ c_{n+1}=c_n+n=\frac{n(n+1)}{2} ~ \text{and}\\
                               & \varepsilon_{n+1}=\varepsilon_n+\frac{\varepsilon'_n}{f}+\big[(n-2)\frac{f'}{f^{2}}+\frac{f''}{ff'}\big](c_n+\varepsilon_n).
\end{aligned}
\end{equation*}
Next, a classical result of Stirling (\cite{13}, p. 151), says
$$\log \Gamma(z)=\left(z-\frac{1}{2}\right)\log z-z+\frac{1}{2}\log (2\pi)+\int_{0}^{+\infty}\frac{[u]-u+\frac{1}{2}}{u+z}du,$$
we have, by Lebesgue's convergence theorem, that
$$f(z)=\frac{\Gamma'}{\Gamma}(z)=\log z-\frac{1}{2z}-\int_{0}^{+\infty}\frac{[u]-u+\frac{1}{2}}{(u+z)^{2}}du,$$
so for $n=1, 2, ...$, one deduces inductively
\begin{equation*}
\begin{aligned}
f^{(n)}(z)=(-1)^{n-1}\left\{\frac{(n-1)!}{z^{n}}+\frac{n!}{2z^{n+1}}+(n+1)!\int_{0}^{+\infty}\frac{[u]-u+\frac{1}{2}}{(u+z)^{n+2}}du \right\}.
\end{aligned}
\end{equation*}
It is thus quite straightforward to verify
\begin{equation}\label{eq1.5}
\begin{aligned}
f(z)=\log z+o(1)=\log z(1+o(1))
\end{aligned}
\end{equation}
and
\begin{equation*}
\begin{aligned}
f^{(n)}(z)=\frac{(-1)^{n-1}(n-1)!}{z^{n}}(1+o(1))
\end{aligned}
\end{equation*}
for $n=1, 2, ...$, and hence
\begin{equation}\label{eq1.6*}
\begin{aligned}
\frac{f'}{f^{2}}(z)=\frac{1}{z(\log z)^{2}}(1+o(1)) ~~\text{and}~~\frac{f''}{ff'}(z)=-\frac{1}{z\log z}(1+o(1)),
\end{aligned}
\end{equation}
uniformly for all $z\in \mathbb{C}\setminus\{z:|\arg z-\pi|\leq\delta\}$ for some $\delta>0$, where $o(1)$ stands for a quantity that goes to $0$ as $|z|\rightarrow +\infty$.

Now, recall $c_n=\frac{n(n-1)}{2}$ and $\varepsilon_1=\varepsilon_2=0$. It can be deduced from (\ref{eq1.6*}) that, as $|z|\rightarrow +\infty$,
\begin{equation*}
\begin{aligned}
\varepsilon_{3}&=\frac{f''}{ff'}=-\frac{1}{z\log z}(1+o(1)),\\
\varepsilon_{4}&=\varepsilon_3+\frac{\varepsilon'_3}{f}+\big(\frac{f'}{f^{2}}+\frac{f''}{ff'}\big)(c_3+\varepsilon_3)=-\frac{4}{z\log z}(1+o(1)),\\
\varepsilon_{5}&=\varepsilon_4+\frac{\varepsilon'_4}{f}+\big(2\frac{f'}{f^{2}}+\frac{f''}{ff'}\big)(c_4+\varepsilon_4)=-\frac{10}{z\log z}(1+o(1)),\\
\vdots\\
\varepsilon_{n}&=-\frac{n(n-1)(n-2)}{6}\frac{1}{z\log z}(1+o(1))~\text{upon assumption, and thus}\\
\varepsilon_{n+1}&=\varepsilon_n+\frac{\varepsilon'_n}{f}+\big[(n-2)\frac{f'}{f^{2}}+\frac{f''}{ff'}\big](c_n+\varepsilon_n)\\
                 &=-\big[\frac{n(n-1)(n-2)}{6}+c_n\big]\frac{1}{z\log z}(1+o(1))\\
                 &=-\big[\frac{n(n-1)(n-2)}{6}+\frac{n(n-1)}{2}\big]\frac{1}{z\log z}(1+o(1))\\
                 &=-\frac{n(n^{2}-1)}{6}\frac{1}{z\log z}(1+o(1)).
\end{aligned}
\end{equation*}
The above formulas illustrate that for $n=1, 2, ...$,
\begin{equation}\label{eq1.6**}
\begin{aligned}
\varepsilon_{n}(z)=o(1) ~~ \text{as} ~~ |z|\rightarrow +\infty.
\end{aligned}
\end{equation}
Noting $\frac{\Gamma^{(n)}}{\Gamma}=f^{n}\big[1+\frac{f'}{f^2}(c_n+\varepsilon_n)\big]$,$\frac{\Gamma^{(l)}}{\Gamma}=f^{l}\big[1+\frac{f'}{f^2}(c_l+\varepsilon_l)\big]$ by using (\ref{eq1.2}) and (\ref{eq1.4}), we get
\begin{equation}\label{eq1.7}
\begin{aligned}
&P_{p_0}\left(\zeta, \zeta', ..., \zeta^{(m)}; 1, \frac{\Gamma^{(n)}}{\Gamma}, \frac{\Gamma^{(l)}}{\Gamma}\right)\\
=&\sum_{q=0}^{M_{p_0}}\sum_{\lambda\in\Lambda_{p_0}\cap\Lambda_{q}^{\ast}}a_{\lambda}(\zeta, \zeta', ..., \zeta^{(m)})\left(\frac{\Gamma^{(n)}}{\Gamma}\right)^{\lambda_{n}}\left(\frac{\Gamma^{(l)}}{\Gamma}\right)^{\lambda_{l}}\\
=&\sum_{q=0}^{M_{p_0}}\sum_{\lambda\in\Lambda_{p_0}\cap\Lambda_{q}^{\ast}}a_{\lambda}(\zeta, \zeta', ..., \zeta^{(m)})f^{n\lambda_{n}}\big[1+\frac{f'}{f^2}(c_n+\varepsilon_n)\big]^{\lambda_{n}}f^{l\lambda_{l}}\big[1+\frac{f'}{f^2}(c_l+\varepsilon_l)\big]^{\lambda_{l}}\\
=&\sum_{q=0}^{M_{p_0}}f^{q}\sum_{\lambda\in\Lambda_{p_0}\cap\Lambda_{q}^{\ast}}a_{\lambda}(\zeta, \zeta', ..., \zeta^{(m)})\big[1+\frac{f'}{f^{2}}(c_{n}+\varepsilon_{n})\big]^{\lambda_{n}}\big[1+\frac{f'}{f^{2}}(c_{l}+\varepsilon_{l})\big]^{\lambda_{l}}.
\end{aligned}
\end{equation}
Then, we define
$$\Lambda_{j}^{\ast\ast}:=\{\lambda\in\Lambda: |\lambda|^{\ast\ast}=j ~\text{with}~|\lambda|^{\ast\ast}:=\lambda_{n}+\lambda_{l}\},$$
and assume that the largest possible $\lambda_{n}+\lambda_{l}$ is $N_{p_0}$. For any $\lambda=(\lambda_{0},\lambda_{n},\lambda_{l})\in\Lambda_{p_0}\cap\Lambda_{q}^{\ast}\cap\Lambda_{j}^{\ast\ast}$, we have
\begin{equation}\label{eq1.7.1}
\left\{
\begin{aligned}
    \lambda_{0}+\lambda_{n}+\lambda_{l}=p_0, \\
    n\lambda_{n}+l\lambda_{l}=q, \\
    \lambda_{n}+\lambda_{l}=j,
\end{aligned}
\right.
\end{equation}
and denote the matrix of coefficients of system \eqref{eq1.7.1} by
\begin{equation*}
\begin{aligned}
B=\left(\begin{array}{cccc}
     1&    1&   1 \\
     0&    n&   l\\
     0&    1&   1\\
\end{array}\right).
\end{aligned}
\end{equation*}
Thus, we have $\det(B)=n-l\neq 0$. By Cramer's Rule, we know $\lambda=(\lambda_{0}, \lambda_{n},\lambda_{l})\in\Lambda_{p_0}\cap\Lambda_{q}^{\ast}\cap\Lambda_{j}^{\ast\ast}$ is uniquely determined by any fixed $p_0, q, j$, and vice versa. So, we can rewrite (\ref{eq1.7}) as
\begin{equation}\label{eq1.8}
\begin{aligned}
&P_{p_0}\bigg(\zeta, \zeta', ..., \zeta^{(m)}; 1, \frac{\Gamma^{(n)}}{\Gamma}, \frac{\Gamma^{(l)}}{\Gamma}\bigg)(z)\\
=&\sum_{q=0}^{M_{p_0}}f^{q}(z)\sum_{j=0}^{N_{p_0}}a_{\lambda\in\Lambda_{p_0}\cap\Lambda_{q}^{\ast}\cap\Lambda_{j}^{\ast\ast}}(\vec{\zeta}) \big[1+H(c_{n}+\varepsilon_{n})\big]^{\lambda_{n}}\big[1+H(c_{l}+\varepsilon_{l})\big]^{\lambda_{l}}(z)\\
=&\sum_{q=0}^{M_{p_0}}f^{q}(z)\sum_{t=0}^{N_{p_0}}\big[b_{q,t}(\vec{\zeta})H^{t}\big](z)\\
=&f^{M_{p_0}}(z)\big[b_{M_{p_0},0}(\vec{\zeta})+b_{M_{p_0},1}(\vec{\zeta})H+\cdots+b_{M_{p_0},N_{p_0}}(\vec{\zeta})H^{N_{p_0}}\big](z)\\
 &+f^{M_{p_0}-1}(z)\big[b_{M_{p_0}-1,0}(\vec{\zeta})+b_{M_{p_0}-1,1}(\vec{\zeta})H+\cdots+b_{M_{p_0}-1,N_{p_0}}(\vec{\zeta})H^{N_{p_0}}\big](z)\\
 &+\cdots+f(z)\big[b_{1,0}(\vec{\zeta})+b_{1,1}(\vec{\zeta})H+\cdots+b_{1,N_{p_0}}(\vec{\zeta})H^{N_{p_0}}\big](z)\\
 &+\big[b_{0,0}(\vec{\zeta})+b_{0,1}(\vec{\zeta})H+\cdots+b_{0,N_{p_0}}(\vec{\zeta})H^{N_{p_0}}\big](z),
\end{aligned}
\end{equation}
with $H=\frac{f'}{f^{2}}$, $(\vec{\zeta})$ being the abbreviation for the vector function $(\zeta, \zeta', ..., \zeta^{(m)})$, and for some $p_0, q, j$, a term with $|\lambda|=p_0, |\lambda|^{\ast}=q, |\lambda|^{\ast\ast}=j$ may not appear in (\ref{eq1.7}), if so, we simply regard the coefficient $a_{\lambda\in\Lambda_{p_0}\cap\Lambda_{q}^{\ast}\cap\Lambda_{j}^{\ast\ast}}(u_{0}, u_{1}, ..., u_{m})$ to be identically zero.

Here, for fixed $p_0, q$, set $u:=(u_0, u_1,..., u_m)$, the polynomials $a_{\lambda\in\Lambda_{p_0}\cap\Lambda_{q}^{\ast}\cap\Lambda_{j}^{\ast\ast}}(u)$, $b_{q, t}(u)$ satisfy the relations as following:
\begin{equation}\label{eq1.8.1}
\begin{aligned}
b_{q,N_{p_0}}(u)=&a_{\lambda\in\Lambda_{p_0}\cap\Lambda_{q}^{\ast}\cap\Lambda_{N_{p_0}}^{\ast\ast}}(u)(c_{n}+\varepsilon_{n})^{\lambda_{n}} (c_{l}+\varepsilon_{l})^{\lambda_{l}}, \\
b_{q,N_{p_0}-1}(u)=&a_{\lambda\in\Lambda_{p_0}\cap\Lambda_{q}^{\ast}\cap\Lambda_{N_{p_0}-1}^{\ast\ast}}(u)(c_{n}+\varepsilon_{n})^{\lambda_{n}}(c_{l}+\varepsilon_{l})^{\lambda_{l}}+a_{\lambda\in\Lambda_{p_0}\cap\Lambda_{q}^{\ast}\cap\Lambda_{N_{p_0}}^{\ast\ast}}(u)\\
          &[\lambda_{l}(c_{n}+\varepsilon_{n})^{\lambda_{n}-1}(c_{l}+\varepsilon_{l})^{\lambda_{l}}+(c_{n}+\varepsilon_{n})^{\lambda_{n}} \lambda_{n}(c_{n}+\varepsilon_{n})^{\lambda_{n}-1}],\\
\vdots\\
b_{q,1}(u)=&\sum_{j=1}^{N_{p_0}}a_{\lambda\in\Lambda_{p_0}\cap\Lambda_{q}^{\ast}\cap\Lambda_{j}^{\ast\ast}}(u)[\lambda_{n}(c_{n}+\varepsilon_{n})+\lambda_{l}(c_{l}+\varepsilon_{l})],\\
b_{q,0}(u)=&\sum_{j=0}^{N_{p_0}}a_{\lambda\in\Lambda_{p_0}\cap\Lambda_{q}^{\ast}\cap\Lambda_{j}^{\ast\ast}}(u).
\end{aligned}
\end{equation}
Recall (\ref{eq1.4}), suppose $b_{q_{0},t_{0}}(u_{0}, u_{1}, ..., u_{m})$ is the first non-zero term in the ordered sequence as following:
$$b_{M_{p_0},0}(u), b_{M_{p_0}-1,0}(u), ..., b_{1,0}(u), b_{0,0}(u),$$
$$b_{M_{p_0},1}(u), b_{M_{p_0}-1,1}(u), ..., b_{1,1}(u), b_{0,1}(u),$$
$$\ddots$$
$$b_{M_{p_0},N_{p_0}-1}(u), b_{M_{p_0}-1,N_{p_0}-1}(u), ..., b_{1,N_{p_0}-1}(u), b_{0,N_{p_0}-1}(u).$$
$$b_{M_{p_0},N_{p_0}}(u), b_{M_{p_0}-1,N_{p_0}}(u), ..., b_{0,N_{p_0}}(u), b_{0,N_{p_0}}(u).$$

Since $b_{q_{0},t_{0}}(u_0, u_1,..., u_m)\not\equiv 0$, we can find an $\varepsilon_0>0$ and a (sufficiently small) subset $\mathbf{\Omega}$ of $\mathbb{C}^{m+1}$ such that,
\begin{equation*}
\begin{aligned}
|b_{q_{0},t_{0}}(u_0, u_1,..., u_m)|\geq\varepsilon_0,
\end{aligned}
\end{equation*}
uniformly for all $u:=(u_0, u_1,..., u_m)\in \mathbf{\Omega}\subsetneq \mathbb{C}^{m+1}$.
Furthermore, since $a_{\lambda\in\Lambda_{p_0}\cap\Lambda_{q}^{\ast}\cap\Lambda_{j}^{\ast\ast}}$ is a polynomial of $u_{0}, u_{1}, ..., u_{m}$, in view of the finiteness of indices, we can find a constant $C_0> 1$ such that
\begin{equation*}
\begin{aligned} |a_{\lambda\in\Lambda_{p_0}\cap\Lambda_{q}^{\ast}\cap\Lambda_{j}^{\ast\ast}}(u_0, u_1,..., u_m)|\leq C_0
\end{aligned}
\end{equation*}
uniformly for all $u:=(u_0, u_1,..., u_m)\in \mathbf{\Omega}\subsetneq \mathbb{C}^{m+1}$
and for all $0\leq p_0 \leq L$, $0\leq q \leq M_{p_0}$, $0\leq j \leq N_{p_0}$.

Then, using Lemma \ref{lem2.1}, there exists a sequence of real numbers $\{y_{k}\}_{k=1}^{+\infty}$ with $|y_{k}|\rightarrow +\infty$ such that $\gamma(y_k)\in \mathbf{\Omega}\subsetneq \mathbb{C}^{m+1}$ when $x=\frac{3}{4}$. It follows from (\ref{eq1.6**}) and (\ref{eq1.8.1}) that there exists a constant $C_1>1$ such that,
for $z_{k}:=\frac{3}{4}+i y_{k}\in \mathbb{C}$, one has
\begin{equation}\label{eq2.1}
\begin{aligned}
|b_{q_{0},t_{0}}(\zeta, \zeta',..., \zeta^{(m)})(z_{k})|\geq\varepsilon_0  \ \ and \ \ |b_{q,t}(\zeta, \zeta',..., \zeta^{(m)})(z_{k})|\leq C_1,
\end{aligned}
\end{equation}
uniformly for all large $k$ and for all $q, t$ with $0\leq q\leq M_{p_0}$, $0\leq t\leq N_{p_0}$. In view of (\ref{eq1.5}) and (\ref{eq1.6*}), we have
$$f^{q}(z_{k})b_{q,t}(\zeta, \zeta',..., \zeta^{(m)})(z_{k})\left(\frac{f'(z_{k})}{f^2(z_{k})}\right)^{t}$$
is equal to
$$b_{q,t}(\zeta, \zeta',..., \zeta^{(m)})(z_{k})\frac{(\log z_{k})^{q-2t}}{(z_{k})^{t}}(1+o(1))$$
when $k\rightarrow +\infty$, where the indices here either satisfy $t=t_{0}$ with $0\leq q\leq q_{0}$ or satisfy $t_{0}<t\leq N_{p_0}$ with $0\leq q\leq M_{p_0}$. As a result, the term
$$f^{q_{0}}(z_{k})b_{q_{0},t_{0}}(\zeta(z_{k}), \zeta'(z_{k})..., \zeta^{(m)}(z_{k}))\left(\frac{f'(z_{k})}{f^2(z_{k})}\right)^{t_{0}},$$
among all possible terms in (\ref{eq1.8}), dominates in growth when $k\rightarrow +\infty$. In fact, for sufficiently large $k$, if $t=t_0$ with $0\leq q< q_{0}$, one derives,
$$\frac{|\log z_{k}|^{q}}{|z_{k}|^{t_0}|\log z_{k}|^{2t_0}}\ll \frac{|\log z_{k}|^{q_0}}{|z_{k}|^{t_0}|\log z_{k}|^{2t_0}},$$
while if $t_{0}<t\leq N_{p_0}$ with $0\leq q\leq M_{p_0}$, one has
$$|\log z_{k}|^{q}\leq |\log z_{k}|^{M_{p_0}+q_0},$$
$$|z_{k}|^{t}|\log z_{k}|^{2t}\gg |z_{k}|^{t_0}|\log z_{k}|^{M_{p_0}+2t_0},$$
and then one also derives,
$$\frac{|\log z_{k}|^{q}}{|z_{k}|^{t}|\log z_{k}|^{2t}}\ll \frac{|\log z_{k}|^{q_0}}{|z_{k}|^{t_0}|\log z_{k}|^{2t_0}}.$$
We then derive from (\ref{eq1.8}) and (\ref{eq2.1}), as $k\rightarrow \infty$, that
\begin{equation}\label{eq2.2}
\begin{aligned}
\left|P_{p_0}\left(\zeta, \zeta', ..., \zeta^{(m)}; 1, \frac{\Gamma^{(n)}}{\Gamma}, \frac{\Gamma^{(l)}}{\Gamma}\right)(z_{k})\right|\geq \frac{\varepsilon_0}{2} \left|\frac{(\log z_{k})^{q_{0}-2t_{0}}}{(z_{k})^{t_{0}}}\right|.
\end{aligned}
\end{equation}

If $p_0=L$, then, by the definition of $p_0$ and ($\ref{eq2.2}$), it yields that
\begin{equation*}
\begin{aligned}
&P(\zeta, \zeta', ..., \zeta^{(m)}; \Gamma, \Gamma^{(n)}, \Gamma^{(l)})(z_{k})\\
=&P_{L}(\zeta, \zeta', ..., \zeta^{(m)}; \Gamma, \Gamma^{(n)}, \Gamma^{(l)})(z_{k})\\
=&\Gamma^{L}(z_{k})P_{L}\left(\zeta, \zeta', ..., \zeta^{(m)}; 1, \frac{\Gamma^{(n)}}{\Gamma}, \frac{\Gamma^{(l)}}{\Gamma}\right)(z_{k})\neq0,
\end{aligned}
\end{equation*}
for sufficiently large $k$. This contradicts the hypothesis (\ref{eq1.1}). Claim \ref{cla 3.1} is proved in this case.

If $p_0<L$, by virtue of another result of Stirling (\cite{13}, p.151) which says
\begin{equation}\label{eq2.2.1}
\begin{aligned}
\left|\Gamma\left(\frac{3}{4}+i y \right)\right|=e^{\frac{-\pi|y|}{2}} |y|^{\frac{1}{4}} \sqrt{2\pi}(1+o(1)), ~ \text{as} ~  y\rightarrow +\infty,
\end{aligned}
\end{equation}
and for all $t>0$,
$f^{q}(z_{k})\left(\frac{f'(z_{k})}{f^{2}(z_{k})}\right)^{t}\rightarrow 0$  as $k\rightarrow \infty$,
one easily observes, in view of (\ref{eq1.3}), (\ref{eq1.4}), (\ref{eq2.2}),(\ref{eq2.2.1}) and (\ref{eq2.1}), that there is a constant $C>0$ such that, for all large $q$,
\begin{equation*}
\begin{aligned}
&\left|\frac{P(\zeta, \zeta', ..., \zeta^{(m)}; \Gamma, \Gamma^{(n)}, \Gamma^{(l)})(z_{k})}{\Gamma^{L}(z_{k})}\right|\\
   =&\left|\sum_{s=p_0}^{L}\frac{1}{\Gamma^{L-s}(z_{k})}P_{s}\left(\zeta, \zeta', ..., \zeta^{(m)}; 1, \frac{\Gamma^{(n)}}{\Gamma}, \frac{\Gamma^{(l)}}{\Gamma}\right)(z_{k}) \right|\\
\geq& e^{\frac{(L-p_0)\pi|y_{k}|}{2}}|y_{k}|^{\frac{1}{4}(p_0-L)}(\sqrt{2\pi})^{p_0-L}\eta \frac{|\log z_{k}|^{q_{0}-2t_{0}}}{|z_{k}|^{t_{0}}}  \\
    &-C_2 e^{\frac{(L-p_0-1)\pi|y_{k}|}{2}}|y_{k}|^{\frac{1}{4}(p_0+1-L)}(\sqrt{2\pi})^{(p_0+1-L)}|\log z_{k}|^{C_3}\rightarrow +\infty,
\end{aligned}
\end{equation*}
as $k\rightarrow \infty$ for some constants $C_2, C_3> 0$ depending on $C_1, L$. Hence,
$$P(\zeta, \zeta', ..., \zeta^{(m)}; \Gamma, \Gamma^{(n)}, \Gamma^{(l)})(z_{k})\neq 0$$
for sufficiently large $z_{k}$, which again contradicts the hypothesis (\ref{eq1.1}). Thus, the proof of Claim \ref{cla 3.1} is completed.

\begin{claim}\label{cla 3.2}
For each $0\leq p\leq L$, when
\begin{equation}\label{eq2.3}
\begin{aligned}
P_{p}(\zeta, \zeta', ..., \zeta^{(m)}; \Gamma, \Gamma^{(n)}, \Gamma^{(l)})(z)\equiv 0
\end{aligned}
\end{equation}
for $z\in \mathbb{C}$, then the polynomial $P_{p}(u_{0}, u_{1}, ..., u_{m}; v_{0}, v_{n}, v_{l})$ vanishes identically.
\end{claim}
\emph{Proof.}
When $p=0$, by definition, $P_{0}(u_{0}, u_{1}, ..., u_{m}; v_{0}, v_{n}, v_{l})$ is a polynomial in $u_{0}, u_{1}, ..., u_{m}$ only; so, from the solution \cite{11} to the question posed by Hilbert, $P_0(\zeta, \zeta', ..., \zeta^{(m)})(z)\equiv 0$ leads to $P_0(u_{0}, u_{1}, ..., u_{m})\equiv 0.$

Henceforth, assume $p>0$. For simplicity, we write $p=p_0$ and use the expression (\ref{eq1.8}) and its associated notations. Next, we prove that $b_{q,t}(u_{0}, u_{1}, ..., u_{m})\equiv 0$ for all $0\leq q\leq M_{p_0}$, $0\leq t\leq N_{p_0}$. To this end, we first prove that each of
$$b_{M_{p_0},0}(u_{0}, u_{1}, ..., u_{m}), b_{M_{p_0}-1,0}(u_{0}, u_{1}, ..., u_{m}), ..., b_{0,0}(u_{0}, u_{1}, ..., u_{m}) $$
must be identically equal to zero. Let's start from $b_{M_{p_0},0}(u_{0}, u_{1}, ..., u_{m})$,
\begin{equation*}
\begin{aligned}
b_{M_{p_0},0}(u_{0}, u_{1}, ..., u_{m})=\sum_{j=0}^{N_{p_0}}a_{\lambda\in\Lambda_{p_0}\cap\Lambda_{M_{p_0}}^{\ast}\cap\Lambda_{j}^{\ast\ast}}(u_{0}, u_{1}, ..., u_{m}),
\end{aligned}
\end{equation*}
a polynomial of $u_{0}, u_{1}, ..., u_{m}$, and assume that it does not vanish identically. Then, following what we have done in Claim \ref{cla 3.1}, one has $(\ref{eq2.1})$ (or its analogue for this newly chosen $P_0$).
Among all the terms of $P_{p_0}\big(\zeta, \zeta', ..., \zeta^{(m)}; 1, \frac{\Gamma^{(n)}}{\Gamma}, \frac{\Gamma^{(l)}}{\Gamma}\big)(z_{k})$ described as in (\ref{eq1.8}), the term
$$f^{M_{p_{0}}}(z_{k})b_{M_{p_{0}},0}(\zeta, \zeta', ..., \zeta^{(m)})(z_{k})\sim b_{M_{p_{0}},0}(\zeta, \zeta', ..., \zeta^{(m)})(z_{k})(\log z_{k})^{M_{p_{0}}},$$
dominates in growth for large $k$, since for all $t>0$, as $k\rightarrow \infty$,
$f^{q}(z_{k})\left(\frac{f'(z_{k})}{f^{2}(z_{k})}\right)^{t}\rightarrow 0$. Thus, analogous to (\ref{eq2.2}), we deduce from (\ref{eq1.8}) and (\ref{eq2.1}) that
\begin{equation*}
\begin{aligned}
 &P_{p_0}(\zeta, \zeta', ..., \zeta^{(m)}; \Gamma, \Gamma^{(n)}, \Gamma^{(l)})(z_k)\\
=&\Gamma^{p_0}(z_{k})P_{p_0}\left(\zeta, \zeta', ..., \zeta^{(m)}; 1, \frac{\Gamma^{(n)}}{\Gamma}, \frac{\Gamma^{(l)}}{\Gamma}\right)(z_{k})\neq0
\end{aligned}
\end{equation*}
for sufficiently large $k$, which contradicts with the hypothesis (\ref{eq2.3}). Therefore,
$b_{M_{p_{0}},0}(u_{0}, u_{1}, ..., u_{m})\equiv 0.$ The next term is $b_{M_{p_{0}}-1,0}(u_{0}, u_{1}, ..., u_{m})$ with $$f^{M_{p_{0}}-1}(z_{k})b_{M_{p_{0}}-1,0}(\zeta, \zeta', ..., \zeta^{(m)})(z_{k})\sim b_{M_{p_{0}}-1,0}(\zeta, \zeta', ..., \zeta^{(m)})(z_{k})(\log z_{k})^{M_{p_{0}}-1},$$
so that one can derive $b_{M_{p_{0}}-1,0}(u_{0}, u_{1}, ..., u_{m})\equiv 0$ in exactly the same manner; repeating this process, we have $b_{q,0}(u_{0}, u_{1}, ..., u_{m})\equiv 0$ for each $0\leq q\leq M_{p_{0}}$. Next, after the elimination of $\frac{f'}{f^2}$, one can perform the preceding procedure again for
$$b_{M_{p_{0}},1}(u_{0}, u_{1}, ..., u_{m}), b_{M_{p_{0}}-1,1}(u_{0}, u_{1}, ..., u_{m}), ..., b_{0,1}(u_{0}, u_{1}, ..., u_{m})$$
and obtain that $b_{q,1}(u_{0}, u_{1}, ..., u_{m})\equiv 0$ for each $0\leq q \leq M_{p_{0}}$. Continuing like this, we get $b_{q,t}(u_{0}, u_{1}, ..., u_{m})\equiv 0$ for all $0\leq q \leq M_{p_{0}}$, $0\leq t \leq N_{p_{0}}$. Thus, it can deduce from $(\ref{eq1.8.1})$ that all the coefficients $a_{\lambda\in\Lambda_{p_0}\cap\Lambda_{q}^{\ast}\cap\Lambda_{j}^{\ast\ast}}(u_{0}, u_{1}, ..., u_{m})$ in $P_{p_0}(u_{0}, u_{1}, ..., u_{m}; v_{0}, v_{n}, v_{l})$ are identically zeros. Therefore,
$$P_{p_0}(u_{0}, u_{1}, ..., u_{m}; v_{0}, v_{n}, v_{l})\equiv 0,$$
and Claim \ref{cla 3.2} is proved completely.

It follows from $(\ref{eq1.2})$ and $(\ref{eq1.3})$ that the proof of Theorem~\ref{M-thm} is a straightforward consequence of Claim \ref{cla 3.1} and Claim \ref{cla 3.2}, so that $(\ref{eq1.1})$ indeed leads to $P(u_{0}, u_{1}, ..., u_{m}; v_{0}, v_{n}, v_{l})\equiv 0$.

\vspace{3mm}

School of Mathematical sciences, Peking University, Beijing, 100871,
P.R. China   e-mail: qiongyanwang@math.pku.edu.com

\end{document}